\documentclass[12pt]{amsart}
\usepackage{amsmath,latexsym,amssymb}

\input xypic
\newcommand{\rar}{\rightarrow}
\newcommand{\lar}{\longrightarrow}

\newtheorem{Theorem}{Theorem}[section]
\newtheorem{theorem}{Theorem}[section]

\newtheorem{lemma}[Theorem]{Lemma}

\newtheorem{corollary}[Theorem]{Corollary}

\newtheorem{proposition}[Theorem]{Proposition}

\newtheorem{conjecture}[Theorem]{Conjecture}

\newtheorem{remark}[Theorem]{Remark}

\newtheorem{example}[Theorem]{Example}

\newtheorem{definition}[Theorem]{Definition}

\def\Ass{\mbox{\rm Ass}}
\def\demo{\noindent{\bf Proof. }}
\def\QED{\hfill$\Box$}
\def\Deg{\mbox{\rm Deg}}
\def\hdeg{\mbox{\rm hdeg}}
\def\Hom{\mbox{\rm Hom}}
\def\Ext{\mbox{\rm Ext}}
\def\Tor{\mbox{\rm Tor}}
\def\depth{\mbox{\rm depth }}
\newcommand{\Rees}{\mbox{${\mathcal R}$}}
\def\gr{\mbox{\rm gr}}

\def\XX{{\bf X}}
\def\xx{{\bf x}}
\def\ff{{\bf f}}
\def\AA{{\bf A}}
\def\BB{{\bf B}}
\def\DD{{\bf D}}

\begin{document}

\title[The Chern Coefficients of Local Rings
]{
The Chern Coefficients of Local Rings}

\author[Wolmer V. Vasconcelos]
{Wolmer V. Vasconcelos}

\thanks{AMS 2000 {\em Mathematics Subject Classification}.
Primary 13A30; Secondary 13B22, 13H10, 13H15.}

\thanks{The  author gratefully acknowledge partial
support from the NSF}

\address{Department of Mathematics, Rutgers University,
Piscataway, New Jersey 08854} \email{vasconce@math.rutgers.edu}
\urladdr{www.math.rutgers.edu/{\textasciitilde}vasconce}


\begin{abstract}
The Chern numbers of the title are the first 
coefficients (after the multiplicities) of the Hilbert functions of various filtrations of
ideals of a local ring $(R, \mathfrak{m})$. For a Noetherian
(good) filtration $\mathcal{A}$ of $\mathfrak{m}$-primary ideals, the
positivity and bounds for $e_1(\mathcal{A})$ are well-studied if $R$
is Cohen-Macaulay, or more broadly, if $R$ is a Buchsbaum ring or
mild generalizations thereof. 
 For arbitrary geometric local domains,
 we introduce techniques based on the
theory of maximal Cohen-Macaulay modules and of extended multiplicity
functions to establish  the meaning of the positivity of
$e_1(\mathcal{A})$,
 and to derive lower and upper bounds for
$e_1(\mathcal{A})$. 
\end{abstract}
\maketitle

\begin{center}

{\em Dedicated to Professor Melvin Hochster on the occasion of his
65th birthday}

\end{center}

\bigskip

\section{Introduction} 
Let $(R, \mathfrak{m})$  be a Noetherian
local ring of dimension $d>0$, and let $I$ be an
$\mathfrak{m}$-primary ideal. One of our goals is to study the set of
$I$-good filtrations of $R$. More concretely,  we will consider the
set of multiplicative, decreasing  filtrations of $R$ 
 ideals, $\mathcal{A}=\{ I_n,
 I_0=R,  I_{n+1}=I I_n, n\gg 0 \}$, integral over the $I$-adic
 filtration, 
conveniently coded in 
  the corresponding
Rees algebra and its associated graded ring
\[ \Rees(\mathcal{A}) = \sum_{n\geq 0} I_nt^n, \quad
\gr_{\mathcal{A}}(R) = \sum_{n\geq 0} I_n/I_{n+1}.
 \]

We will study certain strata of these algebras. For that we will
 focus on the role of the Hilbert polynomial of the Hilbert
function $\lambda(R/I_{n+1})$, $n\gg 0$,
\[ H_{\mathcal{A}}^1(n) =P_{\mathcal{A}}^1(n) = \sum_{i=0}^{d} (-1)^i
e_i(\mathcal{A})
{{n+d-i}\choose{d-i}}, \] particularly of its coefficients
$e_0(\mathcal{A})$ and $e_1(\mathcal{A})$.
Two of our main issues are to establish relationships between the
 coefficients $e_i(\mathcal{A})$, for $i=0,1$, and marginally
 $e_2(\mathcal{A})$.  If $R$ is a
 Cohen-Macaulay ring, there are numerous related developments, noteworthy
 ones being given and discussed in \cite{Elias} and \cite{RV05}.
  The situation is very distinct in the non Cohen-Macaulay
 case. Just to illustrate the issue, suppose $d>2$ and consider a
 comparison between $e_0(I)$ and $e_1(I)$, a subject that has
 received considerable attention. It is often possible to
 pass to a reduction $R\rar S$, with    $\dim S=2$, or even $\dim S=1$, so
 that $e_0(I)= e_0(IS)$ and $e_1(I)=e_1(IS)$. If $R$ is
 Cohen-Macaulay, this is straightforward. However, in general
 the relationship between $e_0(IS)$ and $e_1(IS)$ may involve other
 invariants of $S$, some of which may not be easily traceable all the way to
 $R$.

Our perspective is partly influenced by the interpretation of the
coefficient $e_1$ as a {\em tracking number} (see \cite{DV} for the
original  terminology and 
\cite{ni1}),  that is, as a numerical positional tag of the algebra
$\Rees(\mathcal{A})$ in the set of all such algebras with the same
multiplicity. 
The coefficient $e_1$ under various circumstances is
also called the {\em Chern number or Chern coefficient} of the algebra.

This paper is organized around a list of questions and conjectural
statements about the values of $e_1$ for very general filtrations
associated to the $\mathfrak{m}$--primary ideals of a local
Noetherian ring $(R, \mathfrak{m})$:

\begin{enumerate}
\item (Conjecture 1: the negativity conjecture) For every ideal $J$, generated by a
system of parameters, $e_1(J)<0$ if and only if $R$ is not
Cohen-Macaulay. 
\item (Conjecture 2: the positivity conjecture) For every  $\mathfrak{m}$-primary
ideal $I$, for its integral
closure filtration $\mathcal{A}$
\[ {e}_1(\mathcal{A}) \geq 0.\]
\item (Conjecture 3: the uniformity conjecture)
 For each Noetherian local ring $R$, there exist two functions 
$\ff_l(\cdot)$, $\ff_u(\cdot)$ defined with some extended
multiplicity degree over $R$,
 such that for each  $\mathfrak{m}$-primary ideal
$I$ and any $I$-good filtration $\mathcal{A}$
\[ \ff_l(I) \leq {e}_1(\mathcal{A}) \leq \ff_u(I).\]

\item For any two minimal reductions $J_1$, $J_2$ of  an
$\mathfrak{m}$-primary ideal $I$,
\[ e_1(J_1)=e_1(J_2).\]

\end{enumerate}

For general local rings the conjectures may fail for reasons that
will be illustrated by examples.
We will settle Conjecture 1 for domains that are essentially of
finite type over fields by making use of the existence of special
maximal Cohen-Macaulay modules (Theorem~\ref{e1geoMCM}).
 The lower bound in Conjecture 3 is also
settled for general rings through the use of extended degree
functions (Corollary~\ref{lowere1}).
The upper bound  uses the technique of the Brian\c{c}on-Skoda theorem
on perfect fields (Theorem~\ref{ebounds}). We bring no real
understanding to the last question.

\medskip

We are grateful to Alberto Corso, 
Dan Katz, Claudia Polini, Maria E. Rossi, Rodney Sharp, Bernd
Ulrich and Giuseppe Valla for discussions related to  topics in this
paper.

\section{$e_1$ as a tracking number}

Let $(R, \mathfrak{m})$ be a  Noetherian local domain of
dimension $d$, which is a
quotient of a Gorenstein ring.
For an $\mathfrak{m}$-primary ideal $I$, we are going to consider the
set of all graded subalgebras $\AA$ of the integral closure of
$R[It]$,
\[ R[It] \subset \AA \subset \bar{\AA} = \bar{R[It]}.\]
We will assume that $\bar{\AA}$ is a finite $R[It]$-algebra.
We denote the set of these algebras by $\mathfrak{S}(I)$. If
the algebra 
\[ \AA = \sum I_nt^n\]
comes with a filtration that  is decreasing 
 $I_1\supseteq I_2 \supseteq I_3 \supseteq \cdots $,  it has
an associated graded ring
\[ \gr(\AA) = \sum_{n=0}^{\infty} I_n/I_{n+1}.\]

\begin{proposition} If $\AA$ has the condition $S_2$ of Serre, then
$\{I_n, n\geq 0\}$ is a decreasing filtration.

\end{proposition}

\demo See \cite[Proposition 4.6]{icbook}.

\bigskip

We are going to describe the role of the Hilbert coefficient
$e_1(\cdot)$ in the study of normalization of $R[It]$, following
\cite{ni1} and 
\cite{cnz}. For each $\AA=\sum A_nt^n\in \mathfrak{S}(I)$, we  consider
the Hilbert polynomial (for $n\gg 0)$
\[ \lambda(R/A_{n+1}) = \sum_{i=0}^d (-1)^i e_i(\AA)
{{n+d-i}\choose{d-i}}.
\] The multiplicity $e_0(\AA)$ is constant across the whole set
$\mathfrak{S}(I)$, $e_0(\AA)= e_0(I)$. The next proposition shows the
role of $e_1(\AA)$ in tracking $\AA$ in the set $\mathfrak{S}(I)$.

\begin{proposition} \label{e1tracking}
Let $(R, \mathfrak{m})$ be a normal, Noetherian local domain  which is a
quotient of a Gorenstein ring, and let $I$
be an $\mathfrak{m}$-primary ideal. For algebras $\AA$, $\BB$ of
$\mathfrak{S}(I)$:
\begin{enumerate}
\item If the algebras $\AA$ and $\BB$ satisfy
   $\AA\subset \BB$, then $e_1(\AA)\leq e_1(\BB)$.

\item If $\BB$
is the $S_2$-ification of $\AA$, then $e_1(\AA)=e_1(\BB)$.

\item  If the algebras $\AA$ and $\BB$ satisfy the condition
$S_2$ of Serre and 
   $\AA\subset \BB$, then $e_1(\AA)= e_1(\BB)$ if and only if
   $\AA=\BB$.

\end{enumerate}

\end{proposition}

\demo	The first two assertions  follow directly from the relationship between Krull
dimension and the degree of Hilbert polynomials.
 The exact sequence of graded $R[It]$-modules
\[ 0 \rar \AA \lar \BB \lar \BB/\AA \rar 0 \]
gives that the dimension of $\BB/\AA$ is at most $d-1$. Moreover, 
$\dim \BB/\AA= d-1$ if and only if  its multiplicity is \[\deg(\BB/\AA)=
e_1(\BB)-e_1(\AA)> 0.\]

The last assertion follows because with $\AA$ and $\BB$ satisfying
$S_2$, the quotient $\BB/\AA$ is nonzero, will satisfy $S_1$ and therefore has
Krull dimension $d-1$. \QED

\begin{corollary} Given  a sequence of distinct algebras in
$\mathfrak{S}(I)$, 
\[ \AA_0 \subset \AA_1 \subset \cdots \subset \AA_n= \bar{R[It]},\]
 that satisfy
the condition $S_2$ of Serre, then
\[ n \leq e_1(\bar{R[It]})-e_1(\AA_0)
 \leq e_1(\bar{R[It]})-e_1(I).
\]
\end{corollary}

This highlights the importance of having lower bounds for $e_1(I)$
and upper bounds for $e_1(\bar{R[It]})$. For simplicity we denote the
last coefficient as $\bar{e}_1(I)$. In the Cohen-Macaulay case, for
any parameter ideal $J$, $e_1(J)=0$. Upper bounds for $\bar{e}_1(I)$
were given in \cite{ni1}. For instance, \cite[Theorem 3.2(a),(b)]{ni1}
 shows that if $R$ is a Cohen-Macaulay  algebra of type $t$,
essentially of finite type over a perfect field $k$
and $\delta$ is a non zerodivisor in the Jacobian ideal ${\rm Jac}_k(R)$, then
\[
\overline{e}_1(I) \leq \frac{t}{t + 1}\bigl[(d-1)e_0(I) + e_0(I
+\delta R/\delta R)\bigr]\] \,\, and
\[ \overline{e}_1(I) \leq (d-1)\bigl[e_0(I)
-\lambda(R/\overline{I})\bigr] + e_0(I +\delta R/\delta R).\]

\section{Cohen-Macaulayness and the negativity of $e_1$}

Given the role of the Hilbert coefficient $e_1$ as a tracking number
in the normalization of blowup algebras, it is of interest to know
its signature. 

\medskip

Let $(R, \mathfrak{m})$ be a Noetherian local ring of dimension $d$.
If $R$ is Cohen-Macaulay, then $e_1(J)<0$ for an ideal $J$ generated by a  system of
parameters $x_1, \ldots, x_d$. As a consequence,
for any $\mathfrak{m}$-primary ideal $I$, $e_1(I)
\geq 0$.
If $d=1$, the property $e_1(J)=0$ is characteristic of
Cohen-Macaulayness. For $d\geq 2$, the situation is somewhat
different. Consider the ring $R= k[x,y,z]/(z(x,y,z))$. Then for $T=
H^0_{\mathfrak{m}}(R)$ and $S=k[x,y]= R/T$,  $e_1(R)=e_1(S)=0$. 

\medskip

We are going to argue that the negativity of $e_1(J)$ is an
expression of the lack of Cohen-Macaulayness of $R$ in
numerous classes of rings. To provide a 
framework, we state:

\begin{conjecture}{\rm Let $R$ be a Noetherian local ring that admits
an embedding  into a big Cohen-Macaulay module. Then for a parameter ideal $J$,
$e_1(J)< 0$ if and only if $R$ is not Cohen-Macaulay.
}\end{conjecture}

This places restrictions on $R$, in particular $R$ must be unmixed
and equidimensional. As a matter of fact, we will be concerned
almost exclusively with integral domains  that are essentially of
finite type over a field.

\medskip

We next establish the {\em small} version of the conjecture.

\begin{theorem} \label{e1sMCM} Let $(R,\mathfrak{m})$ be a Noetherian local ring of
dimension $d\geq 2$. Suppose there is an embedding
\[ 0 \rar R \lar E \lar C \rar 0,\] where $E$ is a finitely generated
maximal Cohen-Macaulay $R$-module. If $R$ is not Cohen-Macaulay,
then $e_1(J)<0$ for  any parameter ideal $J$. 
\end{theorem}
\demo We may assume that the residue field of $R$ is infinite. We are
going to argue by induction on $d$. For $d=2$, let $J$ be a parameter
ideal. If $R$ is not Cohen-Macaulay, $\depth C=0$.

Let $J=(x,y)$; we may assume that $x$ is a superficial element for
the purpose of computing $e_1(J)$ and is also superficial relative
to $C$, that is, $x$ is not contained in any associated prime of $C$
distinct from $\mathfrak{m}$.

Tensoring the exact sequence above by $R/(x)$, we get the exact
complex
\[ 0 \rar T = \Tor_1^R(R/(x), C) \lar R/(x) \lar E/xE \lar C/xC \rar
0,\] where $T$ is a nonzero module of finite support. Denote
by $S$ the image of $R'=R/(x)$ in $E/xE $. $S$ is a Cohen-Macaulay ring
of dimension $1$. By the Artin-Rees Theorem, for $n\gg 0$, 
$T\cap (y^n)R'=0$, and therefore from the diagram

\[
\diagram 
0 \rto  & T \cap (y^n)R' \rto\dto & (y^n)R' \rto \dto &(y^n)S
\rto\dto  &0 \\
0 \rto          &  T \rto                 & R' \rto  & S \rto  &0
\enddiagram
\]
the Hilbert polynomial of the ideal $yR'$ is
\[ e_0 n - e_1 = e_0(yS)n + \lambda(T).\]
Thus 
\begin{eqnarray} \label{e1dim1}
 e_1(J) = - \lambda(T)< 0,
\end{eqnarray}
as claimed.

\medskip

Assume now that $d\geq 3$, and let $x$ be a superficial element for
$J$ and for the modules $E$ and $C$. In the exact sequence
\begin{eqnarray} \label{Rmodx}
 0 \rar T = \Tor_1^R(R/(x), C) \lar R'= R/(x) \lar E/xE \lar C/xC \rar
0, 
\end{eqnarray} $T$ is either zero, and we would go on with the induction
procedure, or $T$ is a nonzero module of finite support.

\medskip

Let us point out first an elementary  rule for the calculation of Hilbert
coefficients.
Let $(R, \mathfrak{m})$ be a Noetherian local
ring, and   let $\mathcal{A}=\{I_n, n\geq 0\}$ be a filtration as above.
For a finitely generated $R$-module $M$, denote 
  by $e_i(M)$ the Hilbert coefficients of $M$
for the  filtration $\mathcal{A}M= \{ I_nM, n\geq 0 \}$.

\begin{proposition} \label{abc}
Let 
\[ 0 \rar A \lar B \lar C \rar 0\]
be an exact sequence of finitely generated $R$-modules.
If $r=\dim A< s= \dim B$, then $e_i(B)=e_i(C)$ for $i< s-r$.
\end{proposition}

To continue with the proof,
if in the exact sequence (\ref{Rmodx}) $T\neq 0$, by
Proposition~\ref{abc} we have that $e_1(JR')= e_1(J(R'/T))$, and
the resulting  embedding $R'/T \hookrightarrow E/xE$. By the induction
hypothesis, it suffices to prove that if $R'/T$ is Cohen-Macaulay
then $R'$, and therefore $R$, will be Cohen-Macaulay. This is the
content of \cite[Proposition 2.1]{HuL}. For convenience, we give the
details.

We may assume that $R$ is a complete local ring. Since $R$ is
embedded in a maximal Cohen-Macaulay module, any associated prime of
$R$ is an associated prime of $E$ and therefore it is
equidimensional. Consider the exact sequences
\[ 0 \rar T \rar R' \rar S = R'/T \rar 0, \]

\[ 0 \rar R \stackrel{x}{\lar} R \lar R' \rar 0.\] 
From the first sequence,
taking local cohomology,
\[ 0 \rar H_{\mathfrak{m}}^0(T)=T \lar H_{\mathfrak{m}}^0(R') 
\lar H_{\mathfrak{m}}^0(S) = 0,\]
 since $H_{\mathfrak{m}}^i(T)=0$ for $i>0$ and $S$ is Cohen-Macaulay of dimension
 $\geq 2$; one also has   
$H_{\mathfrak{m}}^1(S)=H_{\mathfrak{m}}^1(R')=0$. From the second sequence, 
since the associated primes of
$R$ have dimension $d$, $H_{\mathfrak{m}}^1(R)$ is a finitely
generated $R$-module. Finally,
 by Nakayama Lemma $H_{\mathfrak{m}}^1(R)=0$, and therefore
 $T=H_{\mathfrak{m}}^0(R')=0$. \QED

\bigskip

We now analyze what is required to extend the proof to big
Cohen-Macaulay cases. We are  going to assume that $R$ is an integral
domain and that $E$ is a big balanced Cohen-Macaulay module (see
\cite[Chapter 8]{BH}, \cite{Sharp81}). 
Embed $R$ into $E$,
\[ 0 \rar R \lar E \lar C \rar 0.\]
The argument above ($d\geq 3)$ will work if in the induction argument we can pick
$x\in J$  superficial for the Hilbert polynomial of $J$, avoiding
the finite set of associated primes of $E$ and all associated primes
of $C$ different from $\mathfrak{m}$. It is this last condition that
is the most troublesome.

There is one case when this can be overcome, to wit, when $R$ is a
complete local ring and $E$ is  countably generated. Indeed, $C$ will
be countably generated and $\Ass(C)$ will be a countable set. The
prime avoidance result of
\cite[Lemma 3]{Burch} allows for the choice of $x$. Let us apply
these ideas in an important case.

\begin{theorem} \label{e1geoMCM} Let $(R,\mathfrak{m})$ be a
Noetherian local integral domain essentially of finite type over a
field. 
 If $R$ is not Cohen-Macaulay, then $e_1(J)<0$
 for  any parameter ideal $J$. 
\end{theorem}

\demo Let $A$ be the integral closure of $R$ and $\widehat{R}$ its
completion. Tensor the embedding $R\subset A$ to obtain
\[ 0 \rar \widehat{R} \lar \widehat{R} \otimes_R A=\widehat{A}. \]
From the properties of pseudo-geometric local rings (\cite[Section
37]{Nagata}), $\widehat{A}$ is a reduced semi-local ring with a
decomposition
\[ \widehat{A} = A_1 \times \cdots \times A_r, \]
where each $A_i$ is a complete local domain, of dimension $\dim R$
and finite over $\widehat{R}$.

For each $A_i$ we make use of \cite[Theorem 3.1]{Griffith76} and
\cite[Proposition 1.4]{Griffith78}
and pick a countably generated big balanced Cohen-Macaulay
$A_i$-module and therefore $\widehat{R}$-module. Collecting the $E_i$
we have an
embedding
\[  \widehat{R} \lar A_1 \times \cdots \times A_r\lar  E= E_1
\oplus \cdots \oplus E_r.\]
As $E$ is a countably generated big balanced Cohen-Macaulay
$\widehat{R}$-module, the argument above shows that if $\widehat{R}$
is not Cohen-Macaulay then $e_1(J\widehat{R})<0$. This suffices to
prove the assertion about $R$. \QED 

\begin{remark}{\em 
There are other classes of local rings admitting big balanced Cohen-Macaulay
modules. A crude way to handle it would be: Let $E$
be such a module and assume it has a set of generators of
cardinality   $s$. Let $X$
be a set of indeterminates of cardinality larger than $s$. 
The local ring $S=R[X]_{\mathfrak{m}[X]}$ is $R$-flat and has the
same depth as $R$. If $E'=S \otimes_RE$ is a big balanced Cohen-Macaulay
$S$-module, with its residue field of cardinality larger than that of
the corresponding module $C$, prime avoidance would again work.
Experts have cautioned that $E'$ may not be balanced. 
}\end{remark}

\begin{example}{\rm We will consider some classes of examples.

\medskip

(i) Let $(R, \mathfrak{m})$ be a regular local ring and let $F$ be a
nonzero (finitely generated) free $R$-module. For any non-free submodule $N$ of $F$,
the idealization (trivial extension) of $R$ by $N$, $S=R\oplus N$ is
a non Cohen-Macaulay local ring. Picking $E=R\oplus F$, 
Theorem~\ref{e1geoMCM}
 implies that for any parameter ideal $J\subset S$, $e_1(J)< 0$.
It is not difficult to give an explicit formula for $e_1(J)$ in this
case.

\medskip

(ii) Let $R= \mathbb{R}+ (x,y)\mathbb{C}[x,y]\subset
\mathbb{C}[x,y]$, for $x,y$ distinct indeterminates. $R$ is not
Cohen-Macaulay but its localization $S$ at the maximal irrelevant ideal
is a Buchsbaum ring. It is easy to verify that
$e_1(x,y)=-1$ and that $e_1(S)=0$ for the $\mathfrak{m}$-adic
filtration of $S$.

Note that $R$ has an isolated singularity. For these rings,
\cite[Theorem 5]{cnz} can be extended (does not require the
Cohen-Macaulay condition), and therefore describes bounds  for
$e_1(\bar{\AA})$ of integral closures. Thus, if $\AA$ is the Rees
algebra of
the parameter ideal $J$, one has
\[ e_1(\bar{\AA})-e_1(J) \leq (d-1+ \lambda(R/L))e_0(J),  \]
where $L$ is the Jacobian of $R$. 

In this example,  one has $d=2$, $\lambda(R/L)=1$,
\[ e_1(\bar{\AA})-e_1(J) \leq 2 e_0(J).  \]

\medskip

(iii) Let $k$ be a field of characteristic zero and let
$f=x^3+y^3+z^3$ be a polynomial of $k[x,y,z]$. Set $A=k[x,y,z]/(f)$
and let $R$ be the Rees algebra of the maximal irrelevant ideal
$\mathfrak{m}$ of
$A$. Using the Jacobian criterion,  $R$ is
normal. Because the reduction number of $\mathfrak{m}$ is $2$, $R$ is
not Cohen-Macaulay.
 Furthermore, it is easy to verify that $R$ is not contained in any
 Cohen-Macaulay domain that is finite over $R$. 
Let
 $S=R_{\mathcal{M}}$, where $\mathcal{M}$ is the irrelevant maximal
 ideal of $R$. 
The first superficial element (in the reduction to dimension two)
can be chosen to be prime. Now one takes the integral closure of $S$,  which will be a
maximal Cohen-Macaulay module.

The argument extends to geometric  domains in any characteristic
 if $\depth R=d-1$.

}\end{example}

\begin{remark}{\rm Uniform lower bounds for $e_1$ are rare but still
exist in special cases. For example, if $R$ is a generalized Cohen-Macaulay
ring, then according to \cite[Theorem 5.4]{GoNi03},
\[ e_1(J)\geq - \sum_{i=1}^{d-1}{d-2\choose
i-1}\lambda(H^i_{\mathfrak{m}}(R)),\]
with equality if $R$ is Buchsbaum.

\medskip

It should be observed that uniform lower bounds may not always exist.
For instance, if $A=k[x,y,z]$, and $R $ the idealization of $(x,y)$,
then for the ideal $J=(x,y,z^n)$, $e_1(J)=-n$.

\medskip

The Koszul homology modules $H_i(J)$ of $J$ is a first place where to look for
bounds for $e_1(J)$. We recall \cite[Theorem 4.6.10]{BH}, that the
multiplicity of $J$ is given by the formula
\[ e_0(J)= \lambda(R/J) - \sum_{i=1}^d(-1)^{i-1} h_i(J),
\] where $h_i(J)$ is the length of $H_i(J)$. The summation term is
non-negative and only vanishes if $R$ is Cohen-Macaulay. Unfortunately
it does not gives a bounds for $e_1(J)$. There is a formula involving
these terms in the special case when $J$ is generated by a
$d$-sequence. Then the corresponding {\em approximation complex} is
acyclic, and the Hilbert-Poincar\'e series of $J$ (\cite[Corollary
4.6]{HSV3}) is
\[ \frac{\sum_{i=0}^d (-1)^i h_i(J)t^i}{(1-t)^d},
\] and therefore
\[ e_1(J) =  \sum_{i=1}^d (-1)^i i h_i(J).
\]

Later  we shall prove the existence of lower bounds more
generally, by making use of extended degree functions.
}\end{remark}

\section{Bounds on $\overline{e}_1(I)$ via the Brian\c{c}on--Skoda
number}

We discuss the role of Brian\c{c}on--Skoda type theorems (see
\cite{AHu}, \cite{LipmanSathaye}) in determining some
relationships between the coefficients $e_0(I)$ and
$\overline{e}_1(I)$. We follow the treatment given in \cite[Theorem
3.1]{ni1} and \cite[Theorem 5]{cnz}, but formulated for the non Cohen-Macaulay case.
 We are going to provide a short
proof along the lines of \cite{LipmanSathaye} for the special case
we need: $\mathfrak{m}$--primary ideals in a local 
ring. The general case is treated by Hochster and Huneke in
\cite[1.5.5 and 4.1.5]{HH}. Let $k$ be a perfect field, let $R$ be
a reduced and equidimensional  $k$--algebra essentially of finite
type, and assume that $R$ is affine with $d={\rm dim} \, R$ or
$(R, \mathfrak{m})$ is local with $d={\rm dim}\, R + {\rm trdeg}_k
(R/\mathfrak{m})$. Recall that the {\it Jacobian ideal}   ${\rm Jac}_k(R)$
of $R$ is defined as the $d$--th Fitting ideal of the module of
differentials $\Omega_k(R)$ -- it can be computed explicitly from
a presentation of the algebra. By varying Noether normalizations
one deduces from \cite[Theorem 2]{LipmanSathaye} that the Jacobian
ideal ${\rm Jac}_k(R)$ is contained in the conductor $R \colon
\overline{R}$ of $R$ (see also \cite{Noether}, \cite[3.1]{AB} and
\cite[2.1]{H}); here $\overline{R}$ denotes the integral closure
of $R$ in its total ring of fractions.

\begin{theorem}\label{GBS}
Let $k$ be a perfect field, let $R$ be a reduced local
 $k$--algebra essentially of finite type with the property $S_2$ of
 Serre, and let
$I$ be an  ideal with a minimal reduction generated by $g$ elements.
Denote by $\DD = \sum_{n\geq 0} D_nt^n$ the $S_2$-ification of $R[
It]$.
 Then for every
integer $n$,
\[
{\rm Jac}_k(R) \, \overline{I^{n+g -1}} \subset D_n.
\]
\end{theorem}
\demo The proof is lifted from \cite[3.1]{ni1}, with the modification
required by the use of $\DD$ at the end.

We may assume that $k$ is infinite. Then, passing to a
minimal reduction, we may suppose that $I$ is generated by  $g$
generators. Let $S$ be a finitely generated
$k$--subalgebra of $R$ so that $R=S_{\mathfrak{p}}$ for some
$\mathfrak{p} \in {\rm
Spec}(S)$, and write $S=k[x_1, \ldots, x_e] = k[X_1, \ldots,
X_e]/{\mathfrak a}$ with ${\mathfrak a}=(h_1, \ldots, h_t)$ an ideal of
height $c$. Notice that $S$ is reduced and equidimensional. Let
$K= (f_1, \ldots, f_g)$ be an $S$--ideal with $K_{\mathfrak{p}}=I$, and
consider the extended Rees ring $B=S[Kt, t^{-1}]$. Now $B$ is a
reduced and equidimensional affine $k$--algebra of dimension
$e-c+1$.

Let $\varphi \colon k[X_1, \ldots, X_e, T_1, \ldots, T_g, U]
\twoheadrightarrow B$ be the $k$--epimorphism mapping $X_i$ to
$x_i$, $T_i$ to $f_it$ and $U$ to $t^{-1}$. Its kernel has height
$c+g$ and contains the ideal ${\mathfrak b}$ generated by $\{ h_i,
T_jU-f_j | 1 \leq i \leq t, 1 \leq j \leq g \}$. Consider the
Jacobian matrix of these generators,
\[
\Theta = \left(
\begin{array}{c|cccc}
\displaystyle\frac{\partial h_i}{\partial X_j} & &   \!\!\!\! 0 & & \\
\hline
& U & & & T_1 \\
* & & \ddots & & \vdots \\
& & & U & T_g
\end{array}
\right).
\]
Notice that $ I_{c+g}(\Theta)\supset I_c(\left(\frac{\partial
h_i}{\partial X_j}\right))U^{g-1}(T_1, \ldots, T_g)$. Applying
$\varphi$ we obtain $ {\rm Jac}_k(B) \supset I_{c+g}(\Theta) B
\supset {\rm Jac}_k(S) K t^{-g+2}$. Thus ${\rm Jac}_k(S) K
t^{-g+2}$ is contained in the conductor of $B$. Localizing at
$\mathfrak{p}$
we see that ${\rm Jac}_k(R) I t^{-g+2}$ is in the conductor of the
extended Rees ring $R[It, t^{-1}]$. Hence for every $n$, ${\rm
Jac}_k(R) \, I \, \overline{I^{n+g-1}} \subset I^{n+1}$, which
yields
\[
{\rm Jac}_k(R) \, \overline{I^{n+g-1}} \subset I^{n+1} \colon I 
 \subset D_{n+1}\colon I = D_n,
\]
as $\gr(\DD)_{+}$ has positive grade. \QED

\bigskip

This result, with an  application of \cite[Theorem~5]{cnz}, gives
the following estimation.

\begin{corollary}
Let $k$ be a perfect field, let $(R, \mathfrak{m})$ be a normal, reduced local
 $k$--algebra essentially of finite type of dimension $d$ and let
$I$ be an  $\mathfrak{m}$-primary ideal. If the Jacobian ideal $L$ of
$R$ is  $\mathfrak{m}$-primary, then for any minimal reduction $J$ of
$I$,
\[ \bar{e}_1(I) - e_1(J) \leq (d+ \lambda(R/L)-1) e_0(I).\]
Moreover, if $L\neq R$ and $\bar{I}\neq \mathfrak{m}$, one replaces $-1$ by $-2$.

\end{corollary}

\section{Lower bounds for $\bar{e}_1$} 
Let $(R, \mathfrak{m})$ be a Noetherian local ring of dimension
$d\geq 1$ and let $I$ be an $\mathfrak{m}$-primary ideal. If $R$ is
Cohen-Macaulay, the original lower bound for $e_1(I)$ was provided by 
Narita (\cite{Narita}) and Northcott (\cite{Northcott}),
\[ e_1(I) \geq e_0(I) - \lambda(R/I).\]
It has been improved in several ways (see a detailed discussion in
\cite{RV07}). For non Cohen-Macaulay rings, estimates for $e_1(I)$
are in a state of flux.

\medskip

We are going to experiment with a special class of non Cohen-Macaulay
rings and methods in seeking lower bounds for $\bar{e}_1(I)$. We are
going to assume that $R$ is a normal domain and the minimal reduction
$J$ of $I$ is generated by a $d$-sequence.
 This is the case of normal
Buchsbaum rings, examples of which can be constructed by a machinery
developed in \cite{EvGr} (see  also \cite{StuckradVogel}).

\medskip

Let $\AA= R[Jt]$ and $\BB=\bar{R[Jt]}$. The corresponding Sally
module $S_{\BB/\AA}$ is defined by the exact sequence
\begin{eqnarray}\label{SallyAB}
 0 \rar B_1t \AA \lar \BB_{+} \lar S_{\AA}(\BB) = \bigoplus_{n\geq
2}B_n/B_1J^{n-1} \rar 0. 
\end{eqnarray}

\begin{lemma} Let $(R, \mathfrak{m})$ be an analytically unramified
local ring and let $J$
be an ideal generated by a system of
parameters $x_1, \ldots, x_d$.
Then
\[ R/\overline{J}\otimes_R \gr_J(R) \simeq R/\overline{J}[T_1,
\ldots, T_d].  \]
\end{lemma}

\demo Let $R[T_1, \ldots, T_d]\rar \AA$ be a minimal presentation of
 $\AA=R[Jt]$.  According to
  \cite[Theorems 3.1, 3.6]{Reesred}, the presentation ideal $\mathcal{L}$ has all
  coefficients in
  $\overline{J}$, that is $\mathcal{L}\subset \bar{J}R[T_1, \ldots,
  T_d]$. \QED

\begin{proposition}
 Let $(R, \mathfrak{m})$ be a normal, analytically unramified
local domain and let $J$
be an ideal generated by a system of
parameters $J=(x_1, \ldots, x_d)$ of linear type.
The Sally module $S_J(\BB)$ defined  above is either $0$ or a
module of dimension $d$ and multiplicity
\[ e_0(S_J(\BB))=\deg(S_J(\BB))= \bar{e}_1(I)-e_0(I) -e_1(J)+\lambda(R/\bar{I}).\]
\end{proposition}

\demo By \cite[Corollary 4.6]{HSV3}, since $R$ is normal of dimension
$d\geq 2$, $\depth \gr_J(R)\geq 2$. Now we make use of \cite[Lemma
1.1]{Hu0} (see also \cite[Proposition 3.11]{icbook}), $\depth
R[Jt]\geq 2$ as well. Now using \cite[Theorem 3.53]{icbook}, we have
that
 $R[Jt]$ has the condition $S_2$ of Serre.

Consider the exact sequence
\[ 0 \rar B_1 R[Jt] \lar R[Jt] \lar  R/\overline{J}\otimes_R
\gr_J(R) \rar 0. \]
By the Lemma, the last algebra is a polynomial ring in $d$ variables.
Therefore $B_1R[Jt]$ satisfies the condition $S_2$ of Serre.
Thus in
 the defining
sequence (\ref{SallyAB}) of $S_J(\BB)$, either $S$ vanishes (and $B_n= B_1J^n$ for $n\geq
2$) or $\dim S_J(\BB) = d$. In the last case, the calculation of the Hilbert
function (see \cite[Remark 2.17]{icbook})
 of $S_J(\BB)$ gives the asserted expression for its multiplicity.
\QED

\section{Existence of general bounds}

 We shall now treat
bounds for $\bar{e}_1(I)$ for several classes of geometric local rings.
Suppose $\dim R = d > 0$.  Let
$I$ be an $\mathfrak{m}$-primary ideal and 
let $\mathcal{A}=\{A_n, n\geq 0\}$ be a filtration integral over the
$I$-adic   filtration. We may assume that
$I=J = (x_1, \ldots, x_d)$ is a parameter ideal.
We will consider some reductions on $R\rar R'$ such that
$e_i(\mathcal{A})= e_i(\mathcal{A}')$, $i=0,1$, for 
$\mathcal{A}'=\{A_nR', n\geq 0\}$.

\medskip

 Using superficial sequences of length $d-1$, $\xx=\{x_1,
\ldots, x_{d-1}\}$ is the technique of choice for Cohen-Macaulay
rings. In general, as in the equality (\ref{e1dim1}), one needs more
control over $H^0_{\mathfrak{m}}(R/(\xx))$. Let us examine first the
case of generalized Cohen-Macaulay local rings by examining the
effect of certain reductions. 

\medskip

\begin{enumerate}

\item
  If $d\geq 2$ and
$R'=R/H^0_{\mathfrak{m}}(R)$, then
$e_i(\mathcal{A})=e_i(\mathcal{A}')$, for $i=0,1$, by
Proposition~\ref{abc}.
In addition $H^i_{\mathfrak{m}}(R)= H^i_{\mathfrak{m}}(R') 
$ for $i\geq 1$.

\item Another property
is that if $d\geq 2$ and $x_1$ is a superficial element for
$\mathcal{A}$ (that is $R$-regular if $d=2$),  then preservation will
hold passing to $R_1=R/(x_1)$. 
As for the lengths of the local cohomology modules, in case $x_1$
is $R$-regular, from
\[ 0 \rar R \stackrel{x}{\lar} R \lar R_1\rar 0\]
we have the exact sequence
\[ 0 \rar  H^0_{\mathfrak{m}}(R_1) \rar  H^1_{\mathfrak{m}}(R)  
\rar  H^1_{\mathfrak{m}}(R)  \rar  H^1_{\mathfrak{m}}(R_1)  \rar
H^2_{\mathfrak{m}}(R)  \rar \cdots,
\] that gives
\begin{eqnarray*} 
\lambda(H^0_{\mathfrak{m}}(R_1)) & \leq  & \lambda(
H^1_{\mathfrak{m}}(R)) \\
\lambda(H^1_{\mathfrak{m}}(R_1)) & \leq & \lambda(
H^1_{\mathfrak{m}}(R)) +
 \lambda( H^2_{\mathfrak{m}}(R)) \\
 & \vdots & \\
 \lambda(H^{d-2}_{\mathfrak{m}}(R_1)) & \leq &
\lambda( H^{d-2}_{\mathfrak{m}}(R))    +
\lambda(H^{d-1}_{\mathfrak{m}}(R)).\\
\end{eqnarray*}   

\item Let us combine the two transformations. Let $T =
H^0_{\mathfrak{m}}(R)$, set $R'= R/T$, let $x\in I$ be a superficial
element for $\mathcal{A}R'$ and set $R_1= R'/(x)$. As $x$ is regular
on $R'$,  we have the exact
sequence
\[ 0 \rar T/xT \lar R/(x) \lar R'/xR'\rar 0,\]
and the associated exact sequence
\[ 0 \rar T/xT \lar H^0_{\mathfrak{m}}(R/(x)) \lar
H^0_{\mathfrak{m}}(R'/xR')
 \rar 0,\]
since $H^1(T/xT)=0$. Note that this gives
\[ R/(x)/H^0_{\mathfrak{m}}(R/(x)) \simeq
R'/xR'/H^0_{\mathfrak{m}}(R'/xR'). \]

There are two consequences to this calculation:
\begin{eqnarray*}
\lambda(H^0_{\mathfrak{m}}(R/(x))) & \leq &
\lambda(H^0_{\mathfrak{m}}(R)) +
\lambda(H^0_{\mathfrak{m}}(R'/xR'))\\  
& \leq & \lambda(H^0_{\mathfrak{m}}(R)) +
\lambda(H^1_{\mathfrak{m}}(R))\\  
\lambda(H^i_{\mathfrak{m}}(R/(x))) & \leq &
\lambda(H^i_{\mathfrak{m}}(R)) +
\lambda(H^{i+1}_{\mathfrak{m}}(R)), \quad 1\leq i \leq d-2\\
\end{eqnarray*}  
\end{enumerate}

\begin{proposition} \label{reddim1} Let $(R, \mathfrak{m})$ be a
generalized Cohen-Macaulay local ring of
positive depth, with $I$ and $\mathcal{A}$ as
above. If $d\geq 2$, consider a sequence of  $d-1$ reductions of the
type $R\rar
R/(x)$,
 and denote by $S$ the ring $R/(x_1, \ldots, x_{d-1})$. Then $\dim S = 1$ and 
\[ \lambda(H^0_{\mathfrak{m}}(S)) \leq 
T(R)= \sum_{i=1}^{d-1}
{{d-2}\choose{i-1}}\lambda(H^i_{\mathfrak{m}}(R)).
\] Moreover, if $R$ is a Buchsbaum ring, equality holds.
\end{proposition}

Let us illustrate another elementary but useful kind of reduction. Let
$(R, \mathfrak{m})$ be a Noetherian  local domain of dimension $d\geq
2$ and let $I$ be
an $\mathfrak{m}$-primary ideal. Suppose $R$ has a finite extension $S$ with the condition $S_2$ of
Serre, 
\[ 0 \rar R \lar S \lar C \rar 0.\]

 Consider the polynomial ring $R[x,y]$ and
tensor the sequence by  $R'=R[x,y]_{\mathfrak{m}[x,y]}$.
Then there are $a,b\in I$
such that 
the ideal generated by  
 polynomial $\ff=ax+by$ has the  following type of primary
 decomposition:
\[ (\ff) = P \cap Q,\]
where $P$ is a minimal prime of $\ff$ and $Q$ is
$\mathfrak{m}R'$-primary. In the sequence
\[ 0 \rar R' \lar S' = R'\otimes_R S \lar C'= R'\otimes_R C \rar 0 \]
we can find $a\in I$ which is superficial for $C$, and pick $b\in I$
so that $(a,b)$ has codimension $2$. Now reduce the second sequence
modulo $\ff=ax+by$. Noting that $\ff$ will be superficial for $C'$,
we will have an exact sequence
\[0 \rar T \lar R'/(\ff) \lar S'/\ff S' \lar C'/\ff C' \rar 0,\]
in which $T$ has finite length and $S'/\ff S' $ is an integral
domain since $a,b $ is a regular sequence in $S$. This suffices to
establish the assertion.

\medskip

Finally we consider generic reductions on $R$. Let $\XX$
be $d\times (d-1)$ matrix $\XX=(x_{ij})$ in $d(d-1)$ indeterminates
and let $C$ be the local ring
$R[\XX]_{\mathfrak{m}[\XX]}$. The filtration $\mathcal{A}C$
has the same Hilbert polynomial as $\mathcal{A}$.
If $I =(x_1, \ldots, x_d)$ we now define the ideal
\[ (f_1, \ldots, f_{d-1})=(x_1, \ldots, x_d) \cdot \XX. \]

\begin{proposition} Let $R$ be an analytically unramified,
 generalized Cohen-Macaulay integral
domain and $I$ and $\mathcal{A}$ defined as above. Then
 $S= C/(f_1, \ldots, f_{d-1})$ is a local ring of dimension one
 such that $\lambda(H^0_{\mathfrak{m}}(S)) \leq T(R)$ and 
$S/H^0_{\mathfrak{m}}(S)$ is an integral domain.

\end{proposition}

The next result shows the existence of bounds for $\bar{e}_1(I)$ as
in \cite{ni1}.

\bigskip

We outline the strategy to find bounds for $e_1(\mathcal{A})$.
Suppose $(R, \mathfrak{m})$ is a local domain essentially of finite
type over a field and $I$ is an $\mathfrak{m}$-primary ideal, and  
denote by $\mathcal{A}$ a filtration as
above. We first achieve a reduction 
to a one dimensional ring $R\rar R'$, where $e_0(I) = e_0(IR')$,
$e_1(\mathcal{A}) = e_1(\mathcal{A}R')$, and in the sequence
\[ 0 \rar T = H^0_{\mathfrak{m}}(R') \lar R' \lar S \rar 0, \]
$T$ is a prime ideal.

\medskip

 Since $S$ is a one-dimensional integral domain, we have he following
 result.

\begin{proposition} \label{reddim1a} In these conditions,
\begin{eqnarray*}
e_1(\mathcal{A})= e_1(\mathcal{A}R') & = & e_1(\mathcal{A}S) -
\lambda(T) \\
e_1(\mathcal{A}S) &\leq & 
\bar{e}_1(I) = \lambda(\bar{S}/S),
\end{eqnarray*}
where $\bar{S}$ is the integral closure of $S$.
\end{proposition}

This shows that to find bounds for $e_1(\mathcal{A})$ one needs to
trace back to the original ring $R$ the properties of $T$ and $S$. As
to $T$, this is realized if $R$
 is a generalized Cohen-Macaulay ring for the reductions described in
 Proposition~\ref{reddim1}.

\begin{theorem}[Existence of Bounds] \label{ebounds} Let $(R, \mathfrak{m})$ be a
local integral domain of dimension $d$,
essentially of finite type over a perfect field, and let  $I$ be an
$\mathfrak{m}$-primary ideal. If $R$ is a generalized Cohen-Macaulay
ring and $\delta$ is a nonzero element of the Jacobian ideal of $R$
then for any $I$-good filtration $\mathcal{A}$ as above,
\[ e_1(\mathcal{A}) < (d-1) e_0(I) +e_0((I,\delta)/(\delta))- T.\]
\end{theorem}

\demo It is a consequence of the proof of \cite[Theorem 3.2(a)]{ni1}.
The assertion there is that
 \[\bar{e}_1(I)
\leq \frac{t}{t + 1} \bigl[ (d-1)e_0(I) + e_0((I +\delta R)/\delta 
R) \bigr],\]
where $t$ is the Cohen-Macaulay type of $R$. Here we apply it to the
reduction in Proposition~\ref{reddim1a}
 \[\bar{e}_1(IS)
<   (d-1)e_0(I) + e_0((I +\delta R)/\delta  
R) ,\]
dropping the term involving $t$, over which we lose control in the
reduction process. Key to the conclusion is the fact that 
 the  element
$\delta$ survives the reduction.
 \QED

\section{Extended degrees and lower bounds for $e_1$} 
The derivation of upper bounds for
$e_1(\mathcal{A})$
above required that $R$ be a generalized Cohen-Macaulay ring. Let us
do away with this requirement by working with the variation of the extended degree
function $\hdeg$ (\cite{DGV}, \cite{hdeg}) labelled $\hdeg_I$ (see
\cite{Linh},
\cite[p. 142]{icbook}). The same method will provide lower bounds for
$e_1(I)$.

\subsection*{Cohomological degrees}
Let $(R, \mathfrak{m})$ be a Noetherian local ring (or a standard
graded algebra over an Artinian local ring) of infinite residue field. We denote by
$\mathcal{M}(R)$ the category of finitely generated $R$-module (or
the corresponding category of graded $R$-modules).

\medskip

A general class of these functions was introduced in \cite{DGV} , while
a prototype was defined earlier in \cite{hdeg}. In his thesis
(\cite{Gunston}), T. Gunston carried out a more formal examination
of such functions in order to introduce 
his own construction of a new cohomological degree. One of the points
that must be taken care of is that of an appropriate {\em generic
hyperplane} section. Let us recall the setting.

\begin{definition}{\rm If $(R,\mathfrak{m})$ is a local ring, a {\em notion of
genericity}\index{notion of genericity} on $\mathcal{M}(R)$ is a
function
\[ U: \{\textrm{\rm isomorphism classes of $\mathcal{M}(R)$} \}
\lar \{ \textrm{\rm non-empty subsets of $\mathfrak{m}\setminus
\mathfrak{m}^2$}\}
\] subject to the following conditions for each $A\in
\mathcal{M}(R)$:
 \begin{itemize}
\item[{\rm (i)}]If $f-g\in \mathfrak{m}^2$ then $f\in U(A)$ if and only if
$g\in U(A)$.
\item[{\rm (ii)}] The set $\overline{U(A)}\subset
\mathfrak{m}/\mathfrak{m}^2$ contains a non-empty Zariski-open subset. 
\item[{\rm (iii)}] If $\depth A> 0$ and $f\in U(A)$, then $f$ is regular on
$A$.
\end{itemize}
}\end{definition}

There is a similar definition for graded modules. We shall usually
switch notation, denoting the algebra by $S$.

\medskip

Another extension  is that associated to an
$\mathfrak{m}$-primary ideal $I$ (\cite{Linh}): A notion of
genericity on $\mathcal{R}$ with respect to $I$ is a function
\[ U: \{\textrm{\rm isomorphism classes of $\mathcal{M}(R)$} \}
\lar \{ \textrm{\rm non-empty subsets of $I\setminus
\mathfrak{m}I$}\}
\] subject to the following conditions for each $A\in
\mathcal{M}(R)$:
 \begin{itemize}
\item[{\rm (i)}]If $f-g\in \mathfrak{m}I$ then $f\in U(A)$ if and only if
$g\in U(A)$.
\item[{\rm (ii)}] The set $\overline{U(A)}\subset
I/\mathfrak{m}I$ contains a non-empty Zariski-open subset. 
\item[{\rm (iii)}] If $\depth A> 0$ and $f\in U(A)$, then $f$ is regular on
$A$.
\end{itemize}

\medskip

 Fixing a notion of
genericity $U(\cdot)$ one has the following extension of the
classical multiplicity.

\begin{definition}{\rm A {\em cohomological
degree}, or {\em extended
multiplicity function},
 is a function \label{Degnu}
\[\Deg(\cdot) : {\mathcal M}(R) \mapsto {\mathbb N},\]
that satisfies the following conditions.
\begin{itemize}
\item[\rm {(i)}]  If $L = \Gamma_{\mathfrak m}(M)$ is the submodule of
elements of $M$ that are annihilated by a 
power of the maximal ideal and $\overline{M} = M/L$, then
\begin{eqnarray}\label{hs0}
\Deg(M) = \Deg(\overline{M}) + \lambda(L),
\end{eqnarray}
where
$ \lambda(\cdot)$ is the ordinary length function.
\item[\rm {(ii)}] (Bertini's rule) 
 If $M$ has positive depth, there is $h\in \mathfrak{m} \setminus
 \mathfrak{m}^2$,
  such that
\begin{eqnarray} \label{Bertini}
\Deg(M) \geq \Deg(M/hM). 
\end{eqnarray}
\item[{\rm (iii)}] (The calibration rule) If $M$ is a Cohen-Macaulay
module, then
\begin{eqnarray} \label{calibration}
\Deg(M) = \deg(M),
\end{eqnarray}
where $\deg(M)$ is the ordinary multiplicity of $M$ 
\end{itemize}
}\end{definition}

In the case of a notion of genericity relative to an
$\mathfrak{m}$-primary ideal $I$, $\deg(M)=e(I;M)$,
  the Samuel's multiplicity of $M$ relative to $I$.

\medskip

 The existence of cohomological degrees in arbitrary
dimensions was established in \cite{hdeg}:

\begin{definition}\label{hdegdef}{\rm Let  $M$ be a finitely generated graded module
over the graded
 algebra $A$ and  $S$  a 
 Gorenstein graded algebra mapping onto $A$, with maximal graded ideal
${\mathfrak m}$.  Set $\dim S=r$,
 $\dim M= d$. 
The {\em homological degree}\index{homological degree}\index{hdeg,
the homological degree} of $M$\label{hdegnu}
is the integer
\begin{eqnarray} 
 \hdeg(M) &=& \deg(M) +
  \label{homologicaldegree} \\
&& \sum\limits_{i=r-d+1}^{r} {{d-1}\choose{i-r+d-1}}\cdot
 \hdeg(\mbox{\rm Ext}^i_S(M,S)).\nonumber 
\end{eqnarray}
This expression becomes more compact when $\dim M=\dim S=d>0$:
\begin{eqnarray}
 \hdeg(M) & = &\deg (M) + \label{homologicaldegree2}\\
&&
 \sum\limits_{i=1}^{d} {{d-1}\choose{i-1}}\cdot
 \hdeg(\mbox{\rm Ext}^i_S(M,S)). \nonumber
\end{eqnarray}
}\end{definition}

\begin{remark}{\rm 
Note that this definition morphs easily into an extended degree noted
$\hdeg_I$ where Samuel multiplicities relative to $I$ are used.
The definition of $\hdeg $ can be extended to any Noetherian local ring $S$ by
setting $\hdeg(M)= \hdeg(\widehat{S}\otimes_SM)$. On other occasions,
we may also assume that the residue field of $S$ is infinite, an
assumption that can be realized by replacing $(S, \mathfrak{m})$ by
the local ring $S[X]_{\mathfrak{m}S[X]}$. In fact, if $X$ is any set
of indeterminates, the localization is still a Noetherian ring, so
the residue field can be assumed to have any cardinality, as we shall
assume in the proofs.
}\end{remark}

\subsection*{Specialization and torsion}  One of the uses of extended
degrees is the following. Let $M$ be a module and $\xx=\{ x_1,
\ldots, x_r\}$ be a superficial sequence for the module $M$ relative
to  an extended degree
$\Deg$. How to estimate the length of $H^0_{\mathfrak{m}}(M)$ in
terms of $M$?

Let us consider the case of $r=1$. Let $H=H^0_{\mathfrak{m}}(M)$ 
and write 
\begin{eqnarray} \label{hmm} 0 \rar H \lar M \lar M' \rar 0.
\end{eqnarray}
Reduction modulo $x_1$ gives 
the exact sequence
\begin{eqnarray} \label{hhmm}
 0 \rar H/x_1H \lar M/x_1 M\lar M'/x_1 M'\rar 0.
\end{eqnarray}
From the first sequence we have $\Deg(M)=\Deg(H)+\Deg(M')$, and from
the second 
\[ 
 \Deg(M/x_1M) - \Deg(H/x_1H) = \Deg(M'/x_1M') \leq \Deg(M').\]

Taking local cohomology of the second exact sequence yields the  short exact
sequence
\[ 0 \rar H/x_1H \lar H^0_{\mathfrak{m}}(M/x_1 M) 
\lar H^0_{\mathfrak{m}}(M'/x_1 M')\rar 0,\]
from we have the estimation
\begin{eqnarray*} \Deg( H^0_{\mathfrak{m}}(M/x_1 M) )& =&
 \Deg(H/x_1H) + \Deg( H^0_{\mathfrak{m}}(M'/x_1 M')) \\ 
& \leq &  \Deg(H/x_1H) + \Deg( M'/x_1 M') \\ 
& \leq &  \Deg(H) + \Deg( M') = \Deg(M). \\ 
\end{eqnarray*}

We resume these observations as:

\begin{proposition} \label{specialtor} Let $M$ be a module and let $\{x_1, \ldots, x_r\}$
be a superficial sequence relative to $M$ and $\Deg$. Then
\[\lambda(H^0_{\mathfrak{m}}(M/(x_1, \ldots, x_r)M)) 
\leq \Deg(M).
\]
\end{proposition}

Now we derive a more precise formula using $\hdeg$. It will be of use
later.

\begin{theorem} \label{Degreddim1} Let $M$ be a module  of
dimension  $d\geq 2$ 
and let $\xx=\{x_1, \ldots, x_{d-1}\}$ be a superficial sequence for
$M$ and $\hdeg$. Then
\[
\lambda(H^0_{\mathfrak{m}}(M/(\xx)M)) \leq
\lambda(H^0_{\mathfrak{m}}(M))+ T(M).
\]
\end{theorem}

\demo 
Consider the exact sequence  
\[ 0 \rar H = H^0_{\mathfrak{m}}(M) \lar M \lar M' \rar 0.\]
We have $\Ext_S^i(M,S) = \Ext_S^i(M',S)$ for $d>i\geq 0$, and
therefore $T(M) = T(M')$. On the other hand, reduction mod $\xx$ gives
\begin{eqnarray*}
\lambda(H^0_{\mathfrak{m}}(M/(\xx)M)) &\leq &
\lambda(H^0_{\mathfrak{m}}(M'/(\xx)M')) +
\lambda(H/(\xx)H) \\
&\leq &
\lambda(H^0_{\mathfrak{m}}(M'/(\xx)M')) +
\lambda(H), 
\end{eqnarray*}
which shows that it is enough to prove the assertion for $M'$.

\medskip

If $d>2$, the argument in the main theorem of \cite{hdeg}
can be used  to pass to $M'/x_1M'$.
This reduces all the way to the case $d=2$.
Write  $h=\xx$. The assertion requires that
$\lambda(H^0_{\mathfrak{m}}(M/hM)) \leq \hdeg(\Ext_S^1(M,S)).$
We have the cohomology exact sequence
\[ \Ext^1_S(M,S) \stackrel{h}{\lar} \Ext^1_S(M,S)
\lar \Ext^2_S(M/hM,S)\lar \Ext^0_S(M,S)=0,
\] where  
\[ \lambda(H^0_{\mathfrak{m}}(M/hM)) =\hdeg(\Ext_S^2(M/hM,S)).\]

If $\Ext_S^1(M,S)$ has finite length the assertion is clear.
Otherwise $L=\Ext_S^1(M,S)$ is a module of dimension $1$ over a
discrete valuation domain $V$ with $h$ for its parameter. By the
fundamental theorem for such modules,
\[ V = V^r \oplus (\bigoplus_{j=1}^s V/h^{e_j}V), \]
so that multiplication by $h$ gives
\[ \lambda(L/hL) = r+s \leq r+ \sum_{j=1}^s e_j = \hdeg(L).\]

An alternative argument at this point is to consider the exact
sequence (we may assume $\dim 
S=1$)
\[ 0 \rar L_0 \lar L \stackrel{h}{\lar} L \lar L/hL \rar 0, \]
where both $L_0$ and $L/hL$ have finite length. If
$H$ denotes the image of the multiplication by $h$ on $L$, dualizing
we have 
the short exact sequence
\[ 0 \rar \Hom_S(L,S) \stackrel{h}{\lar} \Hom_S(L,S) \lar 
\Ext_S^1(L/hL, S) \lar \Ext_S^1(L,S),  \]
which shows that 
\[ \lambda(L/hL) \leq \deg(L) + \lambda(L_0) = \hdeg(L), \]
as desired.
\QED

\bigskip

We now employ the extended degree $\hdeg_I$ to derive lower bounds
for $e_1(I)$. 
We begin by  making a crude comparison between $\hdeg(M)$ and $\hdeg_I(M)$.

\begin{proposition}
Let $(R,\mathfrak{m})$ be a Noetherian local ring and let $I$ be an
$\mathfrak{m}$-primary ideal. Suppose $\mathfrak{m}^r\subset I$. If
$M$ is an $R$-module of dimension $d$, then
\[ \hdeg_I(M) \leq r^d\cdot \deg(M) +  r^{d-1}\cdot(\hdeg(M)-\deg(M)).
\]
\end{proposition}

\demo If $r$ is the index of nilpotency of $R/I$, for any $R$-module
$L$ of dimension $s$,
\[ \lambda(L/(\mathfrak{m}^r)^n L)\geq \lambda(L/I^nL).\]
The Hilbert polynomial of $L$ gives
\[  \lambda(L/(\mathfrak{m}^r)^n L)= \deg(M) {\frac{r^s}{s!}}n^s +
\mbox{\rm lower terms}.\] 
We now apply this estimate to the definition of $\hdeg(M)$, taking
into account that its terms are evaluated  at modules of decreasing
dimension. \QED

\begin{theorem} \label{boundtorsion} Let $(R,\mathfrak{m})$ be a
Noetherian local ring and let $I$ be an $\mathfrak{m}$-primary ideal,
and let $M$ be a finitely generated
$R$-module of dimension $d\geq 1$.
Let $\xx=\{x_1, \ldots, x_r\}$ be a superficial sequence in $I$
relative to $M$ and $\hdeg_I$. Then
\[ \hdeg_I(M/(\xx)M) \leq \hdeg_I(M).\] Moreover, if $r<d$ then
\[ H^0(M/(\xx)M)\leq \hdeg_I(M)-e(I;M).\]

\end{theorem}

If we apply this to $R$, passing to $R'=R/(x_1, \ldots, x_{d-1})$, we
have the estimate for $H^0_{\mathfrak{m}}(R')$ so that the formula
(\ref{e1dim1}) can be used: 

\begin{corollary}[Lower Bound for $e_1(I)$] \label{lowere1} Let $(R, \mathfrak{m})$ be a
Noetherian local ring of dimension $d\geq 1$. If $I$ is an
$\mathfrak{m}$-primary ideal, then
\[ e_1(I) \geq -\hdeg_I(R) + e_0(I).\]$($Note that
$\hdeg_I(R)-e_0(I)$ is the Cohen-Macaulay defficiency of $R$ relative
to the degree function $\hdeg_I$.$)$
\end{corollary}

\end{document}